\documentclass{article}
\usepackage{amsfonts}
\usepackage{amssymb}
\usepackage{amsmath}
\usepackage[latin5]{inputenc}
\numberwithin{equation}{section}
\usepackage[mathscr]{eucal}
\usepackage{eufrak}

\begin{document}
\begin{center}
\Large{\textbf{ON DUAL TIMELIKE - SPACELIKE MANNHEIM PARTNER CURVES IN $ID_{1}^{3}$}}\newline

\large{\"{O}zcan BEKTA\c{S}$^{1} $ and S\"{u}leyman \c{S}ENYURT $^{2} $}

\[ \textbf{\normalsize 1 Department of Mathematics, Arts and Science Faculty, Rize University, Rize-Turkey} \]
\[ \textbf{\normalsize 2 Department of Mathematics, Arts and Science Faculty, Ordu University, Ordu-Turkey.}\]
\[ ozcanbektas1986@hotmail.com; \ \ senyurtsuleyman@hotmail.com, ssenyurt@odu.edu.tr \]
\end{center}
\centerline{}
\begin{abstract}
The first aim of this paper is to define the dual timelike - spacelike  Mannheim partner curves in Dual Lorentzian Space $ID_{1}^{3} $, the second aim of this paper is to obtain the relationships between the curvatures and the torsions of the dual timelike - spacelike  Mannheim partner curves with respect to each other and the final aim of this paper is to get the necessary and sufficient conditions for the dual timelike -spacelike Mannheim partner curves in $ID_{1}^{3} $ .
\end{abstract}

\begin{flushleft}
{\footnotesize \textbf{2000 AMS Subject Classification:} 53B30,51M30,53A35,53A04
\newline \textbf{Keywords:} Mannheim curves, dual Lorentzian Space, curvature, torsion. }
\end{flushleft}

\section{INTRODUCTION}

\indent As is well-known, a surface is said to be ``ruled'' if it is generated by moving a straight line continuously in Euclidean space (O'Neill, 1997). Ruled surfaces are one of the simplest objects in geometric modeling. One important fact about ruled surfaces is that they can be generated by straight lines. A practical application of this type surfaces is that they are used in civil engineering and physics (Guan et al., 1997).

 Since building materials such as wood are straight, they can be considered as straight lines. The results is that if engineers are planning to construct something with curvature, they can use a ruled surface since all the lines are straight (Orbay et al., 2009).

 In the differential geometry of a regular curve in the Euclidean 3 - space $IE^{3} $, it is well-known that one of the important problem is the characterization of a regular curve. The curvature functions $k_{1} $ and $k_{2} $ of a reguler curve play an important role to determine the shape and size of the curve (Kuhnel, 1999; Do Carmo and M.P, 1976). For example, If $k_{1} =k_{2} =0$,  the curve is geodesic. If $k_{1} \ne 0\left(constant\right)$and $k_{2} =0$, then  the curve is a circle with radius  ${1\mathord{\left/ {\vphantom {1 k_{1} }} \right. \kern-\nulldelimiterspace} k_{1} } $. . If $k_{1} \ne 0\left(constant\right)$and $k_{2} \ne 0$$\left(constant\right)$, then the curve is a helix in the space.

 Another way to classification and characterization of curves is the relationship between the Frenet vectors of the curves. For example Saint Venant proposed the question whether upon the surfaces generated by the principal normal of a curve, a second curve can exist which has for its principal normal the principal normal of the given curve. This question was answered by Bertrand in 1850; he showed that a necessary and sufficient condition for the existence of such a second curve is that a linear relationship with constant coefficients exists between the first and second curvatures of the given original curve. The pairs of curves of this kind have been called Conjugate Bertrand curves, or more commonly Bertrand Curves. There are many works related with Bertrand curves in the Euclidean space and Minkowski space. Another kind of associated curves are called Mannheim curve and Mannheim partner curve. If there exists a corresponding relationship between the space curves $\alpha $ and $\beta $ such that, at the corresponding points of the curves, principal normal lines of $\alpha $ coincides with the binormal lines of $\beta $, then $\alpha $ is called a Mannheim curve, and $\beta $ Mannheim partner curve of $\alpha $.

 In recent studies, Liu and Wang (2007,2008) are curious about the Mannheim curves in both Euclidean and Minkowski 3-  space and they obtained

\noindent the necessary and sufficient conditions between the curvature and the torsion for a curve to be the Mannheim partner curves. Meanwhile, the detailed discussion concerned with the Mannheim curves can be found in literature (Wang and Liu, 2007; Liu and Wang, 2008; Orbay and et al., 2009; Özkald$\imath$ et al., 2009; Azak, 2009) and references therein.

 Dual numbers had been introduced by W.K. Clifford (1849 - 1879) as a tool for his geometrical investigations. After him E. Study used dual numbers and dual vectors in his research on line geometry and kinematics. He devoted special attention to the representation of oriented lines by dual unit vectors and defined the famous mapping: The set of oriented lines in an Euclidean three -- dimension space $IE^{3} $ is one to one correspondence with the points of a dual space $ID^{3} $ of triples of dual numbers.

 In this paper, we study the dual timelike - spacelike  Mannheim partner curves in dual Lorentzian space $ID_{1}^{3} $.

\noindent

\section{PRELIMINARY}

\indent By a dual number $A$, we mean an ordered pair of the form $\left(a,a^{*} \right)$ for all $a,a^{*} \in IR$. Let the set $IR\times IR$ be denoted as $ID$. Two inner operations and an equality on $ID=\left\{\left(a,a^{*} \right)\left|a,a^{*} \in IR\right. \right\}$ are defined as follows:

\indent $\left(i\right)\oplus :ID\times ID\to ID$, $A\oplus B=\left(a,a^{*} \right)\oplus \left(b,b^{*} \right)=\left(a+b,a^{*} +b^{*} \right)$ is called the addition in $ID$,\\
\indent $\left(ii\right)\odot :ID\times ID\to ID$. $A\odot B=\left(a,a^{*} \right)\odot \left(b,b^{*} \right)=\left(ab,ab^{*} +a^{*} b\right)$is called the multiplication in $ID$,\\
\indent$\left(iii\right)$ $A=B$ iff $a=b$, $a^{*} =b^{*} $.

\indent If the operations of addition, multiplication and equality on $ID=IR\times IR$ with set of real numbers $IR$are defined as above, the set $ID$ is called the dual numbers system and the element $(a,a^{*})$ of $ID$ is called a dual number. In a dual number $A=(a,a^{*})\in ID$, the real number $a$ is called the real part of $A$ and the real number $a^{*}$ is called the dual part of $A$ The dual number $1=(1,0)$ is called the unit element of multiplication operation $ID$ with respect to multiplication and denoted by $\varepsilon$. In accordance with the definition of the operation of multiplication, it can be easily seen that $\varepsilon^2=0$. Also, the dual number $A=(a,a^{*})\in ID$ can be written as $A=a+\varepsilon a^*$. \\
\indent The set $ID=\{A=a+\varepsilon ^*a|a,a^* \in IR\}$ of dual numbers is a commutative ring according to the operations,\\
\indent i) $(a+\varepsilon a^*)+(b+\varepsilon b^*)=(a+b)+\varepsilon(a^*+b^*)$\\
\indent ii)$(a+\varepsilon a^*)(b+\varepsilon b^*)=ab+\varepsilon(ab^*+ba^*)$. \\
\indent The dual number $A=a+\varepsilon a^*$ divided by the dual number $B=b+\varepsilon b^*$ provided $b \ne 0$ can be defined as \\
\indent $\frac{A}{B}=\frac{a+\varepsilon a^*}{b+\varepsilon b^*}=\frac{a}{b}+\varepsilon \frac{a^*b-ab^*}{b^2}.$ \\
\indent Now let us consider the differentiable dual function. If the dual function $f$ expansions the Taylor series then we have \\
\indent $f(a+\varepsilon a^*)=f(a)+\varepsilon a^*f'(a)$\\
\noindent where $f'(a)$ is the derivation of $f$. Thus we can obtain\\
\indent $sin(a+\varepsilon a^*)=sin a+ \varepsilon a^*cos a$ \\
\indent $cos(a+\varepsilon a^*)=cos a- \varepsilon a^*sin a$\\
\indent The set of $ID^{3} =\{\overrightarrow{A}| \ \ \overrightarrow{A}=\overrightarrow{a}+\varepsilon\overrightarrow{a^{*}},\overrightarrow{a},\overrightarrow{a^{*}}\in IR^3\}$ is a module on the ring $ID$. For any $\overrightarrow{A}=\overrightarrow{a}+\varepsilon\overrightarrow{a^{*}}, \overrightarrow{B}=\overrightarrow{b}+\varepsilon\overrightarrow{b^{*}}\in ID^3$, the scalar or inner product and the vector product of $\overrightarrow{A}$ and $\overrightarrow{B}$ are defined by, respectively,

\indent $\langle \overrightarrow{A}, \overrightarrow{B}\rangle=\langle \overrightarrow{a}, \overrightarrow{b}\rangle+\varepsilon(\langle \overrightarrow{a}, \overrightarrow{b^{*}}\rangle+\langle \overrightarrow{a^{*}}, \overrightarrow{b}\rangle)$,\\
\indent $\overrightarrow{A}\wedge \overrightarrow{B} = \overrightarrow{a} \wedge \overrightarrow{b} + \varepsilon(\overrightarrow{a} \wedge \overrightarrow{b^*}+ \overrightarrow{a^*} \wedge \overrightarrow{b}).$\\
\indent If $\overrightarrow{a} \ne 0$, the norm $\|\overrightarrow{A}\| $ of $\overrightarrow{A}=\overrightarrow{a}+\varepsilon\overrightarrow{a^{*}}$ is defined by \\
\indent $\left\| \overrightarrow{A}\right\| =\sqrt{\left|\left\langle \overrightarrow{A},\overrightarrow{A}\right\rangle \right|} =\left\| \overrightarrow{a}\right\| +\varepsilon \frac{\left\langle \overrightarrow{a},\overrightarrow{a}^{*} \right\rangle }{\left\| \overrightarrow{a}\right\| } {\rm \; },\, {\rm \; }\left\| \overrightarrow{a}\right\| \ne 0 .$ \\
\indent A dual vector $\overrightarrow{A}$ with norm $1$ is called a dual unit vector. The set \\
\indent $S^2=\{\overrightarrow{A}=\overrightarrow{a}+\varepsilon\overrightarrow{a^{*}} \in ID^3 | \|\overrightarrow{A}\|=(1,0), \overrightarrow{a},\overrightarrow{a^{*}}\in IR^3\}$ \\
is called the dual unit sphere with the center $\overrightarrow{O}$ in $ID^3$.\\
\indent Let $\alpha(t)=(\alpha_1(t),\alpha_2(t),\alpha_3(t))$ and $\beta(t)=(\beta_1(t),\beta_2(t),\beta_3(t))$ be real valued curves in $IE^3$. Then $\widetilde{\alpha}(t)=\alpha(t)+\varepsilon\alpha^*(t)$ is a curve in $ID^3$ and it is called dual space curve. If the real valued functions $\alpha_i(t)$ and $\alpha_i^*(t)$ are differentiable then the dual space curve $\widetilde{\alpha}(t)$ is differentiable in $ID^3$. The real part $\alpha(t)$ of the dual space curve $\widetilde{\alpha}=\widetilde{\alpha}(t)$ is called indicatrix. The dual arc-length of real dual space curve $\widetilde{\alpha}(t)$ from $t_1$ to $t$ is defined by \\
\indent $\widetilde{s}=\int_{t_1}^t \| \overrightarrow{\widetilde{\alpha'}}(t) \| dt = \int_{t_1}^t \| \overrightarrow{\alpha{'}}(t) \| dt + \varepsilon = {\int_{t_1}^t \langle \overrightarrow{t}, (\overrightarrow{\alpha^*}(t))^{'}\rangle} dt = s + \varepsilon s^* $

\noindent $\overrightarrow{t}$ is unit tangent vector of the indicatrix $\alpha(t)$ which is a real space curve in $IE^3$. From now on we will take the arc length $s$ of $\overrightarrow{\alpha(t)}$ as the parameter instead of $t$ \\
\indent The Lorentzian inner product of dual vectors $\overrightarrow{A}, \overrightarrow{B}\in ID^3$ is defined by \\
\indent $\langle \overrightarrow{A}, \overrightarrow{B}\rangle = \langle \overrightarrow{a}, \overrightarrow{b}\rangle + \varepsilon(\langle \overrightarrow{a}, \overrightarrow{b^*}\rangle + \langle \overrightarrow{a^*}, \overrightarrow{b}\rangle ) $ \\
with the Lorentzian inner product $\overrightarrow{a}=(a_1,a_2,a_3)$ and $\overrightarrow{b}=(b_1,b_2,b_3)\in IR^3$ \\
\indent $\langle \overrightarrow{a}, \overrightarrow{b}\rangle = -a_1b_1 + a_2b_2 + a_3b_3. $ \\
Thus, $ID^3,\langle , \rangle$ is called the dual Lorentzian space and denoted by $ID^3$. We call the elements of $ID^3$ as the dual vectors. For $\overrightarrow{A}\ne \overrightarrow{0}$. the norm $\|\overrightarrow{A}\|$ of $\overrightarrow{A}$ is defined by $\left\| \overrightarrow{A}\right\| =\sqrt{\left|\left\langle \overrightarrow{A},\overrightarrow{A}\right\rangle \right|} $ . The dual vector $\overrightarrow{A}=\overrightarrow{a}+\varepsilon \overrightarrow{a^{*} }$ is called dual spacelike vector if $\left\langle \overrightarrow{A},\overrightarrow{A}\right\rangle >0$ or $\overrightarrow{A}=0$, dual timelike vector if $\left\langle \overrightarrow{A},\overrightarrow{A}\right\rangle <0$ , dual lightlike vector if $\left\langle \overrightarrow{A},\overrightarrow{A}\right\rangle =0$ for $\overrightarrow{A}\ne 0$. The dual Lorentzian cross-product of $\overrightarrow{A}\, ,\, \overrightarrow{{\rm \; }B}\in ID^{3} \, $is defined by \\
\indent $\overrightarrow{A}\wedge \overrightarrow{B}=\overrightarrow{a}\wedge \overrightarrow{b}+\varepsilon \left(\overrightarrow{a}\wedge \overrightarrow{b}^{*} +\overrightarrow{a}^{*} \wedge \overrightarrow{b}\right)$ \\
where $\overrightarrow{a}\wedge \overrightarrow{b}=\left(a_{3} b_{2} -a_{2} b_{3} ,a_{1} b_{3} -a_{3} b_{1} ,a_{1} b_{2} -a_{2} b_{1} \right)$ $\overrightarrow{a},\overrightarrow{b}\in IR^{3} $i s the Lorentzian cross product.

Dual number \textbf{$\Phi =\varphi +\varepsilon \varphi ^{*} $ } is called dual angle between \textbf{ $\overrightarrow{A}\, \, {\rm ve}\, \, \overrightarrow{B}\, \, $}unit dual vectors. Then we was \\
\indent $\sinh \left(\varphi +\varepsilon \varphi ^{*} \right)=\sinh \varphi +\varepsilon \varphi ^{*} \cosh \varphi $ \\
\indent $\cosh \left(\varphi +\varepsilon \varphi ^{*} \right)=\cosh \varphi +\varepsilon \varphi ^{*} \sinh \varphi.$ \\
\indent Let $\left\{T\left(s\right),N\left(s\right),B\left(s\right)\right\}$be the moving Frenet frame along the curve $\widetilde{\alpha }\left(s\right)$. Then $T\left(s\right)$,$N\left(s\right)$ and $B\left(s\right)$ are dual tangent, the dual principal normal and the dual binormal vector of the curve  $\widetilde{\alpha }\left(s\right)$, respectively. Depending on the casual character of the curve $\widetilde{\alpha }$, we have the following dual Frenet formulas:

If $\widetilde{\alpha }$ is a dual timelike curve ;

\begin{equation} \label{GrindEQ__1_1_}
\left(\begin{array}{c} {T'} \\ {N'} \\ {B'} \end{array}\right)=\left(\begin{array}{ccc} {0} & {\kappa } & {0} \\ {\kappa } & {0} & {\tau } \\ {0} & {-\tau } & {0} \end{array}\right)\left(\begin{array}{c} {T} \\ {N} \\ {B} \end{array}\right)
\end{equation} \\
where $\left\langle T,T\right\rangle =-1,\left\langle N,N\right\rangle =\left\langle B,B\right\rangle =1,\, \left\langle T,N\right\rangle =\left\langle N,B\right\rangle =\left\langle T,B\right\rangle =0.$

\noindent We denote by $\left\{V_{1} \left(s\right),V_{2} \left(s\right),V_{3} \left(s\right)\right\}$ the moving Frenet frame along the curve $\widetilde{\beta }\left(s\right)$. Then $V_{1} \left(s\right),V_{2} \left(s\right)$ and $V_{3} \left(s\right)$ are dual tangent, the dual principal normal and the dual binormal vector of the curve  $\widetilde{\beta }\left(s\right)$, respectively. Depending on the casual character of the curve $\widetilde{\beta }$, we have the following dual Frenet -- Serret formulas:

If $\widetilde{\beta }$ is a dual spacelike curve with a dual timelike binormal $V_{3} $;

\begin{equation} \label{GrindEQ__1_2_}
\left(\begin{array}{c} {V_{1} ^{{'} } } \\ {V'_{2} } \\ {V_{3} ^{{'} } } \end{array}\right)=\left(\begin{array}{ccc} {0} & {P} & {0} \\ {-P} & {0} & {Q} \\ {0} & {Q} & {0} \end{array}\right)\left(\begin{array}{c} {V_{1} } \\ {V_{2} } \\ {V_{3} } \end{array}\right)
\end{equation} \\
where $\left\langle T,T\right\rangle =\left\langle N,n\right\rangle =1,\, \left\langle B,B\right\rangle =-1,\, \left\langle T,N\right\rangle =\left\langle N,B\right\rangle =\left\langle T,B\right\rangle =0.$

\noindent If the curves are unit speed curve, then curvature and torsion calculated by,

\begin{equation} \label{GrindEQ__1_3_}
\left\{\begin{array}{l} {\kappa =\left\| T'\right\| ,} \\ {\tau =\left\langle N^{{'} } ,B\right\rangle ,} \\ {} \\ {P=\left\| V_{1} ^{{'} } \right\| ,} \\ {Q=\left\langle V_{2} ^{{'} } ,V_{3} \right\rangle .} \end{array}\right.
\end{equation}

\noindent If the curves are not unit speed curve, then curvature and torsion calculated by,

\begin{equation} \label{GrindEQ__1_4_}
\left\{\begin{array}{l} {\kappa =\frac{\left\| \widetilde{\alpha }^{{'} } {\wedge \widetilde{\alpha }^{''}}\right\|}{{\left\| \widetilde{\alpha }^{'}{} \right\|}^{3}} ,\, \, \, \, \, \, \, \, \, \, \, \, \, \, \, \, \, \, \, \, \, \, \, \, \, \, \, \, \, \, \, \, \, \, \, \tau =\frac{\det \left(\widetilde{\alpha }^{{'} } ,\widetilde{\alpha }^{{'} {'} } ,\widetilde{\alpha }^{{'} {'} {'} } \right)}{{\left\| \widetilde{\alpha }^{{'} } \wedge \widetilde{\alpha }^{{'} {'} } \right\|}^{2}}  ,} \\ {} \\ {P=\frac{\|{\widetilde{\beta}^{'}\wedge{\widetilde{\beta}^{''}\|}}}{\|{\widetilde{\beta}^{'}}\|^{3}}},\, \, \, \, \, \, \, \, \, \, \, \, \, \, \, \, \, \, \, \, \, \, \, \, \, \, \, \, \, \, \, \, \, \, \, Q =\frac{\det \left(\widetilde{\beta }^{{'} } ,\widetilde{\beta }^{{'} {'} } ,\widetilde{\beta }^{{'} {'} {'} } \right)}{{\left\| \widetilde{\beta }^{{'} } \wedge \widetilde{\alpha }^{{'} {'} } \right\|}^{2}}
\end{array}\right.
\end{equation}

\noindent \textbf{Definition 2.1. a) Dual Hyperbolic angle: }Let $\overrightarrow{A}$ and $\overrightarrow{B}$ be dual timelike vectors in $ID_{1}^{3} $. Then the dual angle between $\overrightarrow{A}$ and $\overrightarrow{B}$ is defined by  $\left\langle \overrightarrow{A},\overrightarrow{B}\right\rangle =-\left\| \overrightarrow{A}\right\| \, \, \left\| \overrightarrow{B}\right\| \cosh \Phi $. The dual number $\Phi =\theta +\varepsilon \theta ^{*} $  is called the dual hyberbolic angle.

\noindent \textbf{b) Dual Central angle: }Let $\overrightarrow{A}$ and $\overrightarrow{B}$ be spacelike vectors in$ID_{1}^{3} $ that span a dual timelike vector subspace. Then the dual angle between $\overrightarrow{A}$ and $\overrightarrow{B}$ is defined by  $\left\langle \overrightarrow{A},\overrightarrow{B}\right\rangle =\left\| \overrightarrow{A}\right\| \, \, \left\| \overrightarrow{B}\right\| \cosh \Phi $. The dual number $\Phi =\theta +\varepsilon \theta ^{*} $  is called  the dual central angle.

\noindent \textbf{c) Dual Spacelike angle: }Let $\overrightarrow{A}$ and $\overrightarrow{B}$ be dual spacelike vectors in$ID_{1}^{3} $ that span a dual spacelike vector subspace. Then the dual angle between $\overrightarrow{A}$ and $\overrightarrow{B}$ is defined by  $\left\langle \overrightarrow{A},\overrightarrow{B}\right\rangle =\left\| \overrightarrow{A}\right\| \, \, \left\| \overrightarrow{B}\right\| \cos \Phi $. The dual number $\Phi =\theta +\varepsilon \theta ^{*} $  is called  the dual spacelike angle.

\noindent \textbf{d) Dual Lorentzian timelike angle: }Let\textbf{ }$\overrightarrow{A}$ be a dual spacelike vector and $\overrightarrow{B}$ be a dual timelike vector in $ID_{1}^{3} $. Then the dual angle between $\overrightarrow{A}$ and $\overrightarrow{B}$ is defined by  $\left\langle \overrightarrow{A},\overrightarrow{B}\right\rangle =\left\| \overrightarrow{A}\right\| \, \, \left\| \overrightarrow{B}\right\| \sinh \Phi $. The dual number $\Phi =\theta +\varepsilon \theta ^{*} $  is called the dual Lorentzian timelike angle.

\section{DUAL TIMELIKE - SPACELIKE MANNHEIM \newline PARTNER CURVE IN $ID_{1}^{3}$}

\indent In this section, we define dual timelike - spacelike Mannheim partner curves in  $ID_{1}^{3} $ and we give some characterization for dual timelike - spacelike Mannheim partner curves in the same space. Using these relationships, we will comment again Shell's and Mannheim's theorems.

\noindent \textbf{Definition 3.1. }Let\textbf{ $\widetilde{\alpha }:I\to ID_{1}^{3} $, $\widetilde{\alpha }\left(s\right)=\alpha \left(s\right)+\varepsilon \alpha ^{*} \left(s\right)$ }be a\textbf{ }dual timelike curve\textbf{ }and \newline\textbf{ $\widetilde{\beta }:I\to ID_{1}^{3} $,$\widetilde{\beta }\left(s\right)=\beta \left(s\right)+\varepsilon \beta ^{*} \left(s\right)$ }be\textbf{ }dual spacelike  with timelike binormal. If there exists a corresponding relationship between the dual timelike curve $\widetilde{\alpha }$ and the dual spacelike curve with dual timelike binormal $\widetilde{\beta }$ such that, at the corresponding points of the curves, the dual binormal lines of $\widetilde{\alpha }$ coincides with the dual principal normal lines of $\widetilde{\beta }$, then $\widetilde{\alpha }$ is called a dual timelike Mannheim curve, and $\widetilde{\beta }$ is called a dual Mannheim partner curve of $\widetilde{\alpha }$. The pair $\left\{\widetilde{\alpha },\widetilde{\beta }\right\}$ is said to be dual timelike - spacelike Mannheim pair. Let $\left\{T,N,B\right\}$ be the dual Frenet frame field along $\widetilde{\alpha }=\widetilde{\alpha }\left(s\right)$ and let $\left\{V_{1} ,V_{2} ,V_{3} \right\}$ be the Frenet frame field along $\widetilde{\beta }=\widetilde{\beta }\left(s\right)$. On the way $\Phi =\theta +\varepsilon \theta ^{*} $ is dual angle between $T$ and $V_{1} $ , there is an following equations between the Frenet vectors and their derivative;

\begin{equation} \label{GrindEQ__2_1_}
\left(\begin{array}{c} {V_{1} ^{{'} } } \\ {V_{2} ^{{'} } } \\ {V_{3} ^{{'} } } \end{array}\right)=\left(\begin{array}{ccc} {\sinh \Phi } & {\cosh \Phi } & {0} \\ {0} & {0} & {1} \\ {\cosh \Phi } & {\sinh \Phi } & {0} \end{array}\right)\left(\begin{array}{c} {T} \\ {N} \\ {B} \end{array}\right).
\end{equation} \\
\textbf{Theorem 3.1. }The distance between corresponding dual points of the dual timelike - spacelike Mannheim partner curves in\textbf{ }$ID_{1}^{3} $ is constant.

\noindent \textbf{Proof: }From the definition of dual spacelike Mannheim curve, we can write

\begin{equation} \label{GrindEQ__2_2_}
\tilde{\beta }(s^{*} )=\; \; \tilde{\alpha }(s)+\lambda \left(s\right)B\left(s\right)
\end{equation}

\noindent By taking the derivate of this equation with respect to $s$ and applying the Frenet formulas, we get

\begin{equation} \label{GrindEQ__2_3_}
V_{1} \frac{ds^{*} }{ds} =T-\lambda \tau N+\lambda 'B
\end{equation}
where the superscript $\left('\right)$ denotes the derivative with respect to the arc length parameter s of the dual curve $\tilde{\alpha }(s)$. Since the dual vectors $B$ and $V_{2} $ are linearly, we get

\indent $\left\langle V_{1} \frac{ds^{*} }{ds} ,B\right\rangle =\left\langle T,B\right\rangle -\lambda \tau \left\langle N,B\right\rangle +\lambda '\left\langle B,B\right\rangle$ and $\lambda '=0$ \\
If we take $\lambda =\lambda _{1} +\varepsilon \lambda _{1}^{*} $, we get $\lambda '_{1} =0$ ve $\lambda _{1}^{*'} =0$ . From here, we can write $\lambda _{1} =c_{1}$ and $\lambda_{1}^{*}=c_{2},$ $c_{1},c_{2}=cons.$

\noindent Then we get $\lambda =c_{1} +\varepsilon c_{2} $. On the other hand, from the definition\textbf{ }of distance function between $\tilde{\alpha }(s)$ and $\tilde{\beta }(s)$ we can write \\
\indent $d\left(\tilde{\alpha }(s),\tilde{\beta }(s)\right)=\left\| \tilde{\beta }(s)-\tilde{\alpha }(s)\right\| = \left|\lambda _{1} \right|\mp \varepsilon \lambda _{1}^{*} = \left|c_{1} \right|\mp \varepsilon c_{2}$ \\
This is completed the proof.

\noindent \textbf{Theorem 3.2. }For a dual timelike - spacelike curve\textbf{ }$\widetilde{\alpha }$ in \textbf{ }$ID_{1}^{3} $, there is a dual spacelike curve $\tilde{\beta }$ so that $\left\{\tilde{\alpha },\tilde{\beta }\right\}$ is a dual spacelike Mannheim pair.\\
\noindent \textbf{Proof: }Since\textbf{ }the dual vectors $V_{2} $ and $B$ are linearly dependent, the equation \eqref{GrindEQ__2_2_} can be written as

\begin{equation} \label{GrindEQ__2_4_}
\tilde{\alpha }=\tilde{\beta }-\lambda V_{2}
\end{equation}
Since $\lambda $ is a nonzero constant, there is a dual timelike curve $\tilde{\beta }$ for all values of $\lambda $.

 Now, we can give the following theorem related to curvature and torsion of the dual timelike - spacelike Mannheim partner curves.

\noindent \textbf{Theorem 3.3. }Let $\left\{\tilde{\alpha },\tilde{\beta }\right\}$ be a dual timelike - spacelike Mannheim pair in \textbf{$ID_{1}^{3} $. }If\textbf{ }$\tau $ is dual torsion of \textbf{ }$\tilde{\alpha }$ and $P$ is dual curvature and $Q$ is dual torsion of  $\tilde{\beta }$ , then\textbf{}

\begin{equation} \label{GrindEQ__2_5_}
\tau =-\frac{P}{\lambda Q}
\end{equation}
\textbf{Proof: }By taking the derivate of equation\textbf{ }\eqref{GrindEQ__2_3_}  with respect to $s$ and applying the Frenet formulas, we obtain

\begin{equation} \label{GrindEQ__2_6_}
V_{1} \frac{ds^{*} }{ds} =T-\lambda \tau N
\end{equation}
Let \textbf{$\Phi =\varphi +\varepsilon \varphi ^{*} $ }be\textbf{ }dual angle between\textbf{ } the dual tangent vectors $T$ and $V_{1} $, we can write

\begin{equation} \label{GrindEQ__2_7_}
\left\{\begin{array}{l} {V_{1} =\sinh \Phi \, T+\cosh \Phi \, N} \\ {V_{3} =\cosh \Phi \, T+\sinh \Phi \, N} \end{array}\right.
\end{equation}
From \eqref{GrindEQ__2_6_} and \eqref{GrindEQ__2_7_} , we get

\begin{equation} \label{GrindEQ__2_8_}
\frac{ds^{*} }{ds} =\frac{1}{\sinh \Phi } ,\, \, \, \, \, -\lambda \tau =\cosh \Phi \frac{ds^{*} }{ds}
\end{equation}
By taking the derivate of equation\textbf{ }\eqref{GrindEQ__2_4_}  with respect to $s$ and applying the Frenet formulas, we obtain

\begin{equation} \label{GrindEQ__2_9_}
T=\left(1+\lambda P\right)V_{1} \frac{ds^{*} }{ds} -\lambda QV_{3} \frac{ds^{*} }{ds}
\end{equation}
From equation \eqref{GrindEQ__2_7_} we can write

\begin{equation} \label{GrindEQ__2_10_}
\left\{\begin{array}{l} {T=-\sinh \Phi \, V_{1} +\cosh \Phi \, V_{3} } \\ {N=\cosh \Phi \, V_{1} -\sinh \Phi \, V_{3} } \end{array}\right.
\end{equation}
where $\Phi $ is the dual angle between $T$ and $V_{1} $ at the corresponding points of the dual curves of $\tilde{\alpha }$  and $\tilde{\beta }$ . By taking into consideration equations \eqref{GrindEQ__2_9_} and \eqref{GrindEQ__2_10_}, we get

\begin{equation} \label{GrindEQ__2_11_}
\sinh \Phi =-\left(1+\lambda P\right)\frac{ds^{*} }{ds} ,\, \, \, \, \, \cosh \Phi =-\lambda Q\frac{ds^{*} }{ds}
\end{equation}
Substituting $\frac{ds^{*} }{ds} $  into \eqref{GrindEQ__2_11_} , we get

\begin{equation} \label{GrindEQ__2_12_}
\sinh ^{2} \Phi =-\left(1+\lambda P\right),\, \, \, \, \, \cosh ^{2} \Phi =\lambda ^{2} \tau Q
\end{equation}
From the last equation, we can write \\
\indent $\tau =-\frac{P}{\lambda Q}$ \\
If the last equation is seperated into the dual and real parts, we can obtain

\begin{equation} \label{GrindEQ__2_13_}
\left\{\begin{array}{l} {k_{2} =-\frac{p}{cq} } \\ {k_{2}^{*} =\frac{pq^{*} -p^{*} q}{cq^{2} } } \end{array}\right.
\end{equation}
\textbf{Corollary 3.1. }Let $\left\{\tilde{\alpha },\tilde{\beta }\right\}$ be a dual timelike - spacelike Mannheim pair in \textbf{$ID_{1}^{3} $. }Then, the dual product of torsions\textbf{ }$\tau $ and $Q$ at the corresponding points of the dual spacelike Mannheim partner curves is not constant.

 Namely, Schell's theorem is invalid for the dual timelike - spacelike Mannheim curves. By considering Theorem 3.3 we can give the following results.

\noindent \textbf{Corollary 3.2. }Let $\left\{\tilde{\alpha },\tilde{\beta }\right\}$ be a dual timelike - spacelike Mannheim pair in \textbf{$ID_{1}^{3} $. }Then,\textbf{ }torsions\textbf{ }$\tau $\textit{ }and $Q$ has a negative sign.

\noindent \textbf{Theorem 3.4.} Let $\left\{\tilde{\alpha },\tilde{\beta }\right\}$ be a dual timelike - spacelike Mannheim pair in \textbf{$ID_{1}^{3} $. }Between the curvature and the torsion of the dual spacelike curve $\widetilde{\beta }$ , there is the relationship

\begin{equation} \label{GrindEQ__2_14_}
\mu Q-\lambda P=1
\end{equation}
where $\mu $ and$\lambda $ are nonzero dual numbers.

\noindent \textbf{Proof: }From equation \eqref{GrindEQ__2_11_}, we obtain \\
\indent $\frac{\sinh \Phi }{1+\lambda P} =\frac{\cosh \Phi }{\lambda Q}$, \\
arranging this equation, we get \\
\indent $\tanh \Phi =\frac{1+\lambda P}{\lambda Q}$, \\
and if we choose $\mu =\lambda \tanh \Phi $ for brevity, we see that \\
\indent $\mu Q-\lambda P=1$.\\
\textbf{Theorem 3.5. }Let $\left\{\tilde{\alpha },\tilde{\beta }\right\}$ be a dual timelike - spacelike Mannheim pair in \textbf{$ID_{1}^{3} $. }There are the following equations for the curvatures and the torsions of the curves $\widetilde{\alpha }$ ve  $\widetilde{\beta }$ \\
\indent $i) \kappa = -\frac{d\Phi }{ds},$ \\
\indent $ii)\tau = P\cosh \Phi \frac{ds^{*} }{ds} -Q\sinh \Phi \frac{ds^{*} }{ds},$ \\
\indent $iii)P = \tau \cosh \Phi \frac{ds}{ds^{*}}, $ \\
\indent $iv)Q = \tau \sinh \Phi \frac{ds}{ds^{*}}. $

\textbf{Proof: }$i)$By considering equation \eqref{GrindEQ__2_7_}, we can easily that $\left\langle T,V_{1} \right\rangle =\cos \Phi $. Differentiating of this equality with respect to  \textit{s }by considering equation \eqref{GrindEQ__1_1_} , we have \\
\indent $\left\langle T',V_{1} \right\rangle +\left\langle T,V_{1} ^{{'} } \right\rangle =-\sinh \Phi \frac{d\Phi }{ds}$,\\
from equations  \eqref{GrindEQ__1_1_} and \eqref{GrindEQ__1_2_}, we can write \\
\indent $\left\langle \kappa N,V_{1} \right\rangle +\left\langle T,PV_{2} \frac{ds^{*} }{ds} \right\rangle =-\sinh \Phi \frac{d\Phi }{ds}$, \\
from equations  \eqref{GrindEQ__2_10_}, we get \\
\indent $\kappa =-\frac{d\Phi }{ds}$.\\
If the last equation is seperated into the dual and real part, we can obtain \\
$ii)$ By considering equation \eqref{GrindEQ__2_7_}, we can easily that $\left\langle N,V_{2} \right\rangle =0$. Differentiating of this equality with respect to  \textit{s and }by considering equation \eqref{GrindEQ__1_1_} , we have \\
\indent $\left\langle N',V_{2} \right\rangle +\left\langle N,V_{2} ^{{'} } \frac{ds^{*} }{ds} \right\rangle =0$, \\
From equations  \eqref{GrindEQ__1_1_} and \eqref{GrindEQ__1_2_}, we can write \\
\indent $\left\langle \kappa T+\tau B,V_{2} \right\rangle +\left\langle \cosh \Phi \, V_{1} -\sinh \Phi \, V_{3} ,\left(-PV_{1} +QV_{3} \right)\frac{ds^{*} }{ds} \right\rangle =0$,\\
From equations  \eqref{GrindEQ__2_10_}, we get \\
\indent $\tau =P\cosh \Phi \frac{ds^{*} }{ds} -Q\sinh \Phi \frac{ds^{*} }{ds}$, \\
$iii)$ By considering equation \eqref{GrindEQ__2_7_}, we can easily that $\, \left\langle B,V_{1} \right\rangle =0$. Differentiating of this equality with respect to  \textit{s and }by considering equation \eqref{GrindEQ__1_1_} , we have \\
\indent $\left\langle B',V_{1} \right\rangle +\left\langle B,V_{1} ^{{'} } \frac{ds^{*} }{ds} \right\rangle =0$, \\
From equations  \eqref{GrindEQ__1_1_}, \eqref{GrindEQ__1_2_} and \eqref{GrindEQ__2_10_} we can write \\
\indent $\left\langle -\tau \left(\cosh \Phi \, V_{1} -\sinh \Phi \, V_{3} \right),V_{1} \right\rangle +\left\langle B,PV_{2} \frac{ds^{*} }{ds} \right\rangle =0$,\\
\indent $P=\tau \cosh \Phi \frac{ds}{ds^{*} }$,\\
$iv)$ By considering equation \eqref{GrindEQ__2_7_}, we can easily that $\left\langle B,V_{3} \right\rangle =0$. Differentiating of this equality with respect to  \textit{s }by considering equation \eqref{GrindEQ__1_1_} , we have \\
\indent $\left\langle B',V_{3} \right\rangle +\left\langle B,V_{3} ^{{'} } \frac{ds^{*} }{ds} \right\rangle =0$,\\
From equations  \eqref{GrindEQ__1_1_}, \eqref{GrindEQ__1_2_} and \eqref{GrindEQ__2_10_} we can write \\
\indent $\left\langle -\tau \left(\cosh \Phi \, V_{1} -\sinh \Phi \, V_{3} \right),V_{3} \right\rangle +\left\langle B,QV_{2} \frac{ds^{*} }{ds} \right\rangle =0$,\\
\indent $Q=\tau \sinh \Phi \frac{ds}{ds^{*} }$. \\
\textbf{Corollary 3.3. }Let $\left\{\tilde{\alpha },\tilde{\beta }\right\}$ be a dual timelike - spacelike Mannheim pair in \textbf{$ID_{1}^{3} $. }If the statements of Theorem 3.5 is seperated into the dual and real part, we can obtain \\
\indent $i)\left\{\begin{array}{l} {k_{2} =p\cosh \theta \frac{ds^{*} }{ds} -q\sinh \theta \frac{ds^{*} }{ds} } \\ {k_{2}^{*} =\left(p^{*} \cosh \theta +p\theta ^{*} \sinh \theta \right)\frac{ds^{*} }{ds} -\left(q^{*} \sinh \theta +q\theta ^{*} \cosh \theta \right)\frac{ds^{*} }{ds} } \end{array}\right. $

\indent $ii)\left\{\begin{array}{l} {p=k_{2} \cosh \theta \frac{ds}{ds^{*} } } \\ {p^{*} =\left(k_{2}^{*} \cosh \theta +k_{2} \theta ^{*} \sinh \theta \right)\frac{ds}{ds^{*} } ,} \end{array}\right. $

\indent $iii)\left\{\begin{array}{l} {q=k_{2} \sinh \theta \frac{ds}{ds^{*} } } \\ {q^{*} =\left(k_{2}^{*} \sinh \theta +k_{2} \theta ^{*} \cosh \theta \right)\frac{ds}{ds^{*} } .} \end{array}\right. $ \\
\noindent By considering the statements iii) and iv) of Theorem 2.5 we can give the following results.

\noindent \textbf{Corollary 3.4. }Let $\left\{\tilde{\alpha },\tilde{\beta }\right\}$ be a dual timelike - spacelike Mannheim pair in \textbf{$ID_{1}^{3} $. }Then there exist the following relation between curvature and torsion of\textbf{ }$\widetilde{\beta }$ and torsion of $\widetilde{\alpha }$;

\begin{equation} \label{GrindEQ__2_15_}
P^{2} -Q^{2} =\tau ^{2} \left(\frac{ds}{ds^{*} } \right)^{2}
\end{equation}
\textbf{Theorem 3.6. }A dual timelike space curve in\textbf{ }$ID_{1}^{3} $ is a dual timelike - spacelike Mannheim curve if and only if its curvature $P$ and torsion $Q$ satisfy the formula

\begin{equation} \label{GrindEQ__2_16_}
\lambda \left(Q^{2} -P^{2} \right)=P
\end{equation}
where $\lambda $ is never pure dual constant.

\noindent \textbf{Proof:} By taking the derivate of the statement $\widetilde{\alpha }=\widetilde{\beta }-\lambda V_{2} $ with respect to $s$ and applying the Frenet formulas we obtain \\
\indent $T\frac{ds}{ds^{*} } =V_{1} +\lambda \left(PV_{1} -QV_{3} \right)$, \\
\indent $\kappa N\left(\frac{ds}{ds^{*2} } \right) +T\frac{d^{2} s}{ds^{*2}}=PV_{2} +\lambda \left(P'V_{1} -Q'V_{3} +\left(P^{2} -Q^{2} \right)V_{2} \right)$ \\
Taking the inner product the last equation with $B$, we get \\
\indent $\lambda \left(Q^{2} -P^{2} \right)=P$. \\
If the last equation is seperated into the dual and real part, we can obtain

\begin{equation} \label{GrindEQ__2_17_}
\left\{\begin{array}{l} {p=\lambda \left(q^{2} -p^{2} \right)} \\ {p^{*} =2\lambda \left(qq^{*} -pp^{*} \right)} \end{array}\right.
\end{equation}
where $\lambda =c_{1} +\varepsilon c_{2} $ .

\noindent \textbf{Theorem 3.7. }Let $\left\{\tilde{\alpha },\tilde{\beta }\right\}$ be a dual timelike - spacelike Mannheim partner curves in \textbf{$ID_{1}^{3} $. }Moreover, the dual points\textbf{ }$\widetilde{\alpha }\left(s\right)$, $\widetilde{\beta }\left(s\right)$ be two corresponding dual points of $\left\{\widetilde{\alpha },\widetilde{\beta }\right\}$and $M$ ve $M^{*} $ be the curvature centers at these points, respectively. Then, the ratio  \\

\begin{equation} \label{GrindEQ__2_18_}
\frac{\left\| \widetilde{\beta }\left(s\right)M\right\| }{\left\| \widetilde{\alpha }\left(s\right)M\right\| } :\frac{\left\| \widetilde{\beta }\left(s\right)M^{*} \right\| }{\left\| \widetilde{\alpha }\left(s\right)M^{*} \right\| }=\left(1+\kappa P\right)\left(1+\lambda P\right)\ne constant.
\end{equation}

\noindent \textbf{Proof:  }A circle that lies in the dual osculating plane of the point $\widetilde{\alpha }\left(s\right)$ on the dual timelike curve $\widetilde{\alpha }$ and that has the centre $M=\widetilde{\alpha }\left(s\right)+\frac{1}{\kappa } N$ lying on the dual principal normal $N$ of the point $\widetilde{\alpha }\left(s\right)$ and the radius $\frac{1}{\kappa } $ far from $\widetilde{\alpha }\left(s\right)$, is called dual osculating circle of the dual curve $\widetilde{\alpha }$ in the point $\widetilde{\alpha }\left(s\right)$. Similar definition can be given fort he dual curve $\widetilde{\beta }$ too.

\noindent Then, we can write \\
\indent $\left\| \widetilde{\alpha }\left(s\right)M\right\| =\left\| \frac{1}{\kappa } N\right\| \, \, =\frac{1}{\kappa }$,\\
\indent $\left\| \widetilde{\alpha }\left(s\right)M^{*} \right\| =\left\| \lambda B+\frac{1}{P} V_{2} \right\| =\frac{1}{P} +\lambda$, \\
\indent $\left\| \widetilde{\beta }\left(s\right)M^{*} \right\| =\left\| \frac{1}{P} V_{2} \right\| =\frac{1}{P}$, \\
\indent $\left\| \widetilde{\beta }\left(s\right)M\right\| =\left\| \lambda V_{3} +\frac{1}{\kappa } N\right\| =\frac{1}{\kappa } +\lambda $ \\
Therefore, we obtain \\
\indent $\frac{\left\| \widetilde{\beta }\left(s\right)M\right\| }{\left\| \widetilde{\alpha }\left(s\right)M\right\| } :\frac{\left\| \widetilde{\beta }\left(s\right)M^{*} \right\| }{\left\| \widetilde{\alpha }\left(s\right)M^{*} \right\| } =\left(1+\lambda P\right)\sqrt{1-\lambda ^{2} \kappa ^{2} } \ne cons. $ \\
Thus, we can give the following

\noindent \textbf{Corollary 3.5.} Mannheim's  Theorem is invalid for the dual timelike - spacelike Mannheim partner curve  $\left\{\tilde{\alpha },\tilde{\beta }\right\}$in \textbf{$ID_{1}^{3} $.}

\noindent \textbf{REFERENCES}

\noindent \textbf{[1] }A. Z. Azak, On timelike Mannheim partner curves in $L^{3} $, Sakarya University Faculty of Arts and Science The Journal of Arts and Science, Vol. 11\eqref{GrindEQ__2_}, (2009), 35-45.

\noindent \textbf{[2] }B. O'Neill,  Semi--Riemannian Geometry with Applications to Relativity, Academic Press, New York, 1983.

\noindent \textbf{[3]} B.\textbf{ }O'Neill,~~Elemantary Differential Geometry, 2nd ed. Academic Press, New York, 1997.

\noindent \textbf{[4] }Do Carmo, Manfredo, Differential Geometry of Curves and Surfaces, Pearson Education. New York: Academic Press, 1976.

\noindent \textbf{[5] }H. Liu and F. Wang, Mannheim Partner Curves in 3-space, J. Geom. Vol. 88(1-2), 2008, 120-126.

\noindent \textbf{[6] }H. Liu, F. Wang (2007). Mannheim Partner Curves in 3-space, Procedings of The Eleventh International Workshop on Diff. Geom.11,25-31

\noindent \textbf{[7]}K\textbf{. }Orbay, E. Kasap and İ. Aydemir, Mannheim Offsets of Ruled Surfaces, Mathematical Problems in Engineering, Article Number:160917, 2009.

\noindent \textbf{[8] }K. Orbay, E. Kasap, On Mannheim partner curves in $E^{3} $ , International Journal of Physical Sciences, Vol.4 \eqref{GrindEQ__5_}, 2009, 261-264.

\noindent \textbf{[9]} M. A. Gungor and M. Tosun, A study on dual Mannheim partner curves, International Mathematical Forum 5, no. 45-48, 2010, 2319--2330.

\noindent \textbf{[10] }M. Kazaz, M. Önder, Mannheim Offsets of Timelike Ruled Surfaces in Minkowski 3-space  $IR_{1}^{3} $ , eprint/arXiv:0906.2077v3. 3 1 R .

\noindent \textbf{[11] }R. Blum, A Remarkable Class of Mannheim-Curves, Canad. Math. Bull., Vol. 9\eqref{GrindEQ__2_}, 1966, 223-228.

\noindent \textbf{[12] }S.\textbf{ }Özkaldi, K. İlarslan and Y. Yayli,~~On Mannheim Partner curves in Dual Space,~~Analele Stiintifice ale Universitatii Ovidius Constanta, Seria Matematica, vol XVII, fasc. 2, 2009.

\noindent \textbf{[13] }W.\textbf{ }Kuhnel, Differential Geometry: Curves-Surfaces-Manifolds, Braunschweig, Wiesbaden, 1999.

\noindent \textbf{[14] }Z.\textbf{ }Guan, J. Ling, X. Ping and T. Rongxi (1997). Study and Application of Physics-Based Deformable Curves and surfaces, Computers and Graphics 21: 305-313.

\end{document}